\input amstex
\input amsppt.sty
\magnification=\magstep1
\hsize=32truecc
 \vsize=23.5truecm
\baselineskip=14truept
 \NoBlackBoxes
\def\q{\quad}
\def\qq{\qquad}
\def\mod#1{\ (\text{\rm mod}\ #1)}

\def\t{\text}
\def\qtq#1{\q\t{#1}\q}
\par\q
\def\f{\frac}
\def\e{\equiv}

\def\b{\binom}

\def\o{\omega}

\def\({\left(}
\def\){\right)}
\def\sls#1#2{(\f{#1}{#2})}
 \def\ls#1#2{\big(\f{#1}{#2}\big)}
\def\Ls#1#2{\Big(\f{#1}{#2}\Big)}

\let \pro=\proclaim
\let \endpro=\endproclaim

\topmatter
\title Recurrence formulas with gaps for Bernoulli and
Euler polynomials
\endtitle
\author ZHI-Hong Sun\endauthor
\affil School of Mathematical Sciences, Huaiyin Normal University,
\\ Huaian, Jiangsu 223001, PR China
\\ E-mail: zhihongsun$\@$yahoo.com
\\ Homepage: http://www.hytc.edu.cn/xsjl/szh
\endaffil

 \nologo \NoRunningHeads
\abstract{Let $\{B_n\}$, $\{B_n(x)\}$ and $\{E_n(x)\}$ be the
Bernoulli numbers, Bernoulli polynomials and Euler polynomials,
respectively. In this paper we mainly establish formulas for
$\sum_{6\mid k-3}\binom nkB_{n-k}(x)$, $\sum_{6\mid k}\binom
nkE_{n-k}(x)$ and $\sum_{6\mid k-3}\binom nkm^kB_{n-k}$ in the cases
$m=2,3,4$.
\par\q
\newline MSC: Primary 11B68, Secondary 11B39 \newline Keywords:
Bernoulli number; Bernoulli polynomial; Euler polynomial}
 \endabstract
\endtopmatter
\document
\subheading{1. Introduction}
\par The Bernoulli numbers $\{B_n\}$ and Bernoulli polynomials
$\{B_n(x)\}$ are defined by
$$B_0=1,\ \sum_{k=0}^{n-1}\b nkB_k=0\ (n\ge 2)\qtq{and}
B_n(x)=\sum_{k=0}^n\b nkB_kx^{n-k}\ (n\ge 0).$$  The Euler numbers
$\{E_n\}$ and Euler polynomials $\{E_n(x)\}$ are
 defined by
$$\f{2\t{e}^t}{\t{e}^{2t}+1}=\sum_{n=0}^{\infty}E_n\f{t^n}{n!}\
(|t|<\f{\pi}2)\q \t{and}\q
\f{2\t{e}^{xt}}{\t{e}^t+1}=\sum_{n=0}^{\infty}E_n(x)\f{t^n}{n!}\
(|t|<\pi),\tag 1.1$$ which are equivalent to (see [MOS])
$$E_0=1,\ E_{2n-1}=0,\ \sum_{r=0}^n\b{2n}{2r}E_{2r}=0\ (n\ge 1)$$
and $$E_n(x)+\sum_{r=0}^n\b nrE_r(x)=2x^n\
 (n\ge 0).\tag 1.2$$
 The first few Bernoulli and Euler numbers are shown below:
$$\align&B_0=1,\ B_1=-\f 12,\ B_2=\f 16,\ B_4=-\f 1{30},\ B_6=\f
1{42},\ B_8=-\f 1{30}, B_{10}=\f 5{66},
\\&E_0=1,\ E_2=-1,\ E_4=5,\ E_6=-61,\ E_8=1385,\ E_{10}=-50521.
\endalign$$
 It is well known that ([MOS])
 $$\aligned E_n(x)&=\f 1{2^n}\sum_{r=0}^n\b nr(2x-1)^{n-r}E_r
\\&=\f 2{n+1}\Big(B_{n+1}(x)-2^{n+1}B_{n+1}\Big(\f
x2\Big)\Big)\\&=\f{2^{n+1}}{n+1}\Big(B_{n+1}\Ls{x+1}2-B_{n+1}\Ls
x2\Big).\endaligned\tag
 1.3$$
 In particular,
$$E_n=2^n E_n\Ls 12\qtq{and}E_n(0)=\f{2(1-2^{n+1})B_{n+1}}{n+1}.
\tag 1.4$$
\par In his first paper Ramanujan [R] found some recurrence formulas with gaps for
Bernoulli numbers. In particular, he showed that for
$n=3,5,7,\ldots,$
$$\sum\Sb k=0\\6\mid k-3\endSb^n\b nkB_{n-k}
=\cases -\f n6&\t{if $n\e 1\mod 6$,}
\\\f n3&\t{if $n\e 3,5\mod 6$.}
\endcases\tag 1.5$$
Some generalizations of (1.5) were given by M. Chellali [C] in 1988.
However, the proofs of (1.5) given by Ramanujan and  Chellali are
somewhat complicated. Let $\omega=(-1+\sqrt{-3})/2$. In this paper,
we prove the following generalization of (1.5) for Bernoulli
polynomials in a very simple manner:
$$\sum\Sb
k=0\\6\mid k-3\endSb^n\b nkB_{n-k}(x)=\f
n6\big(x^{n-1}+(x-1)^{n-1}-(x+\omega)^{n-1}-(x+\omega^2)^{n-1}\big).$$
As consequences we give explicit formulas for $\sum\limits\Sb
k=0\\6\mid k-3\endSb^n\b nkm^kB_{n-k}$ in the cases $m=2,3,4$.
\par Let $\Bbb N$ be the set of positive integers, and let $[x]$ be the greatest integer not exceeding $x$. In [L] Lehmer showed that for $n=2,4,6,\ldots$,
$$4E_n=2\big(1+(-3)^{\f n2}\big)-3\sum_{k=1}^{[n/6]}\b
n{6k}2^{6k}E_{n-6k}.\tag 1.6$$ In this paper we prove that for
$n\in\Bbb N$,
$$4E_n(x)+3\sum_{k=1}^{[n/6]}
\b n{6k}E_{n-6k}(x)=x^n+(x-1)^n+(x+\omega)^n+(x+\omega^2)^n.$$
 \subheading{2. Recurrence
formulas with gaps for Bernoulli polynomials}
\par For two numbers $b$ and $c$, let $\{U_n(b,c)\}$ and $\{V_n(b,c)\}$ be
the Lucas sequences given by
$$U_0(b,c)=0,\ U_1(b,c)=1,\ U_{n+1}(b,c)=bU_n(b,c)-cU_{n-1}(b,c)\ (n\ge 1)$$
and
$$V_0(b,c)=2,\ V_1(b,c)=b,\ V_{n+1}(b,c)=bV_n(b,c)-cV_{n-1}(b,c)\ (n\ge 1).$$
 It is well
known that (see [W]) for $b^2-4c\not=0$,
$$U_n(b,c)= \f
1{\sqrt{b^2-4c}}\Big\{\Big(\f{b+\sqrt{b^2-4c}}2\Big)^n
-\Big(\f{b-\sqrt{b^2-4c}}2\Big)^n\Big\}$$ and
$$V_n(b,c)=\Ls{b+\sqrt{b^2-4c}}2^n+\Ls{b-\sqrt{b^2-4c}}2^n.$$
 \pro{Lemma
2.1} For $n\in\Bbb N$ we have
$$\sum_{k=0}^{n-1}\b nkB_k(x)((1+z)^{n-k}-z^{n-k})=n(x+z)^{n-1}.$$
\endpro
Proof. It is known that (see [MOS])
$$\sum_{k=0}^n\b
nkB_k(x)y^{n-k}=B_n(x+y)\qtq{and}B_n(t+1)-B_n(t)=nt^{n-1}.\tag 2.1$$
Thus,
$$\align &\sum_{k=0}^{n-1}\b nkB_k(x)((1+z)^{n-k}-z^{n-k})
\\&=\sum_{k=0}^n\b nkB_k(x)(1+z)^{n-k}-B_n(x)-\sum_{k=0}^n\b nkB_k(x)z^{n-k}+B_n(x)
\\&=B_n(x+z+1)-B_n(x+z)=n(x+z)^{n-1}.
\endalign$$ This proves the lemma.

 \pro{Theorem 2.1} For $n\in\Bbb N$ we have $$\sum\Sb
k=0\\6\mid k-3\endSb^n\b nkB_{n-k}(x)=\f
n6\big(x^{n-1}+(x-1)^{n-1}-(x+\omega)^{n-1}-(x+\omega^2)^{n-1}\big).$$
\endpro
Proof. As $1+\o+\o^2=0$ and $\o^3=1$, using Lemma 2.1 we see that
$$\align&-4\sum\Sb k=0\\n-k\e 3\mod 6\endSb^n\b nkB_k(x)
+2\sum\Sb k=0\\n-k\e \pm 1\mod 6\endSb^n\b nkB_k(x)
\\&
=\sum_{k=0}^{n-1}\b
nkB_k(x)((-1)^{n-k}-1)(\omega^{n-k}+\omega^{2(n-k)})
\\&=\sum_{k=0}^{n-1}\b
nkB_k(x)((-\o^2)^{n-k}-\o^{n-k}+(-\o)^{n-k}-\o^{2(n-k)})
\\&=\sum_{k=0}^{n-1}\b
nkB_k(x)((1+\o)^{n-k}-\o^{n-k}+(1+\o^2)^{n-k}-\o^{2(n-k)}))
\\&=n(x+\o)^{n-1}+n(x+\o^2)^{n-1}.\endalign$$
That is,
$$-6\sum\Sb k=0\\6\mid n-k-3\endSb^n\b nkB_k(x)
+2\sum\Sb k=0\\2\nmid n-k\endSb^n\b nkB_k(x)
=n(x+\o)^{n-1}+n(x+\o^2)^{n-1}.$$ From Lemma 2.1 we have
$$\sum_{k=0}^{n-1}\b
nkB_k(x)=nx^{n-1}\ \t{and}\ \sum_{k=0}^{n-1}\b
nkB_k(x)(-1)^{n-k+1}=n(x-1)^{n-1}.\tag 2.2$$ Thus,
$$2\sum\Sb k=0\\2\nmid n-k\endSb^n\b nkB_k(x)
=\sum_{k=0}^{n-1}\b nkB_k(x)(1-(-1)^{n-k})=nx^{n-1}+n(x-1)^{n-1}.$$
Hence
$$\align 6\sum\Sb k=0\\6\mid k-3\endSb^n\b nkB_{n-k}(x)
&=6\sum\Sb r=0\\6\mid n-r-3\endSb^n\b nrB_r(x)
\\&=nx^{n-1}+n(x-1)^{n-1}-n(x+\o)^{n-1}-n(x+\o^2)^{n-1}.
\endalign$$
This completes the proof.
\par\q
\par As $B_n(0)=B_n$, taking $x=0$ in
Theorem 2.1 we obtain Ramanujan's identity (1.5).

 \pro{Theorem 2.2} For $m\in\Bbb N$ and $n=3,5,7,\ldots$ we have
 $$\align&\sum\Sb
k=0\\6\mid k-3\endSb^n\b nkm^kB_{n-k} \\&=\f
n3\Big\{\sum_{r=1}^{m-1}r^{n-1}-\sum_{r=1}^{[\f
{m-1}2]}\big((r+m\omega)^{n-1}+(r+m\o^2)^{n-1}\big)+\delta(m,n)\Big\},\endalign$$
where $$\delta(m,n)=\cases -\f 12m^{n-1}&\t{if $2\nmid m$ and $3\mid
n-1$,}\\m^{n-1}&\t{if $2\nmid m$ and $3\nmid n-1$,}
\\-\f 12m^{n-1}-(-3)^{\f{n-1}2}\sls m2^{n-1}&\t{if $2\mid m$ and $3\mid
n-1$,}\\m^{n-1}-(-3)^{\f{n-1}2}\sls m2^{n-1}&\t{if $2\mid m$ and
$3\nmid n-1$.}\endcases$$
\endpro
Proof. From Theorem 2.1 we have
$$\align&\sum\Sb
k=0\\6\mid k-3\endSb^n\b nk\sum_{r=0}^{m-1}B_{n-k}\Ls rm
\\&=\f n6
\sum_{r=0}^{m-1}\Big\{\Ls
rm^{n-1}+\Ls{r-m}m^{n-1}-\Ls{r+m\o}m^{n-1}-\Ls
{r+m\o^2}m^{n-1}\Big\}.\endalign$$ By Raabe's theorem (see [MOS]),
$\sum_{r=0}^{m-1}B_{n-k}\sls rm=m^{1-(n-k)}B_{n-k}$. Thus,
$$\align&\sum\Sb
k=0\\6\mid k-3\endSb^n\b nkm^kB_{n-k}\\&=\f n6
\sum_{r=0}^{m-1}\Big\{r^{n-1}+(m-r)^{n-1}-(r+m\o)^{n-1}-(r+m\o^2)^{n-1}\Big\}
\\&=\f n6\Big\{2\sum_{r=1}^{m-1}r^{n-1}+m^{n-1}
-m^{n-1}(\o^{n-1}+\o^{2(n-1)})
\\&\qq-\sum_{r=1}^{m-1}((r+m\o)^{n-1}+(r+m\o^2)^{n-1})\Big\}
.\endalign$$ As $2\mid n-1$ and $1+\o+\o^2=0$, we see that
$$(m-r+m\o)^{n-1}+(m-r+m\o^2)^{n-1}
=(r+m\o)^{n-1}+(r+m\o^2)^{n-1}.$$ Hence
$$\align&\sum\Sb
k=0\\6\mid k-3\endSb^n\b nkm^kB_{n-k}\\&=\f n6
\Big\{m^{n-1}(1-\o^{n-1}-\o^{2(n-1)})-(1+(-1)^m)\Ls
m2^{n-1}(-3)^{\f{n-1}2}
\\&\q+2\sum_{r=1}^{m-1}r^{n-1}-2\sum_{r=1}^{[(m-1)/2]}
\Big((r+m\o)^{n-1}+(r+m\o^2)^{n-1}\Big)\Big\}.
\endalign$$
This yields the result.
\par Putting $m=2$ in Theorem 2.2 we deduce the following result.
 \pro{Corollary 2.1} For $n=3,5,7,\ldots$
we have
$$\sum\Sb k=0\\6\mid k-3\endSb^n\b nk2^kB_{n-k}
=\cases \f n3(1-2^{n-2}-(-3)^{\f{n-1}2})&\t{if $n\e 1\mod 6$,}
\\\f n3(1+2^{n-1}-(-3)^{\f{n-1}2})&\t{if $n\e 3,5\mod 6$.}
\endcases$$
\endpro
\pro{Corollary 2.2} For $n=3,5,7,\ldots$ we have
$$\sum\Sb k=0\\6\mid k-3\endSb^n\b nk3^kB_{n-k}
=\cases \f n3\big(1+2^{n-1}-\f{3^{n-1}}2-V_{n-1}(1,7)\big)&\t{if
$6\mid n-1$,}\\\f n3\big(1+2^{n-1}+3^{n-1}-V_{n-1}(1,7)\big)&\t{if
$6\nmid n-1$.}\endcases$$
 \endpro
Proof. Observe that
$$\align V_{n-1}(1,7)&=\Ls{1+\sqrt{1^2-4\cdot 7}}2^{n-1}+\Ls{1-\sqrt{1^2-4\cdot 7}}2^{n-1}
\\&=\Ls{1+3\sqrt{-3}}2^{n-1}+\Ls{1-3\sqrt{-3}}2^{n-1}
\\&=(1+3\o)^{n-1}+(1+3\o^2)^{n-1}.\endalign$$ Taking $m=3$ in Theorem
2.2 we deduce the result.

 \pro{Corollary 2.3} For $n=3,5,7,\ldots$ we have
$$\aligned&\sum\Sb k=0\\6\mid k-3\endSb^n\b nk4^kB_{n-k}
\\&=\cases \f
n3\big(1+2^{n-1}+3^{n-1}-2^{2n-3}-(-12)^{\f{n-1}2}-V_{n-1}(2,13)
\big)&\t{if $6\mid n-1$,}\\\f
n3\big(1+2^{n-1}+3^{n-1}+4^{n-1}-(-12)^{\f{n-1}2}-V_{n-1}(2,13)\big)&\t{if
$6\nmid n-1$.}\endcases\endaligned$$
\endpro
Proof. Note that
$$\align V_{n-1}(2,13)&=\Ls{2+\sqrt{2^2-4\cdot 13}}2^{n-1}+\Ls{2-\sqrt{2^2-4\cdot 13}}2^{n-1}
\\&=(1+2\sqrt{-3})^{n-1}+(1-2\sqrt{-3})^{n-1}
\\&=(1+4\o)^{n-1}+(1+4\o^2)^{n-1}.\endalign$$ Taking $m=4$ in Theorem
2.2 we deduce the result.

\pro{Theorem 2.3} For $n\in\Bbb N$ we have
$$\sum\Sb k=0\\4\mid k-2\endSb^n\b
nkB_{n-k}(x)(-1)^{\f{k-2}4}2^{n+1-\f
k2}=ni\big((2x-1-i)^{n-1}-(2x-1+i)^{n-1}\big).$$
\endpro
Proof. Using Lemma 2.1 we see that
$$\align&4\sum\Sb k=0\\4\mid n-k-2\endSb^{n-1}\b nkB_k(x)\Ls
i2^{\f{n-k}2}
\\&=\sum\Sb k=0\\2\mid n-k\endSb^{n-1}\b nkB_k(x)\cdot 2\Big(\Ls
i2^{\f{n-k}2}-\Big(-\f i2\Big)^{\f{n-k}2}\Big)
\\&=\sum_{k=0}^{n-1}\b nkB_k(x)(1+(-1)^{n-k})
\Big(\Ls{1+i}2^{n-k}-\Ls{i-1}2^{n-k}\Big)
\\&=\sum_{k=0}^{n-1}\b nkB_k(x)\Big(\Ls{1+i}2^{n-k}-\Ls{i-1}2^{n-k}
-\Ls{1-i}2^{n-k}+\Ls{-i-1}2^{n-k}\Big)
\\&=n\Big(x+\f{i-1}2\Big)^{n-1} -n\Big(x+\f{-i-1}2\Big)^{n-1}.\endalign$$
Thus,
$$\align&\sum\Sb k=0\\4\mid
k-2\endSb^n\b nkB_{n-k}(x)(-1)^{\f{k-2}4}2^{-\f{k-2}2}
\\&=\sum\Sb k=0\\4\mid
k-2\endSb^n\b nkB_{n-k}(x)\Ls i2^{\f{k-2}2}=\sum\Sb r=0\\4\mid
n-r-2\endSb^n\b nrB_r(x)\Ls i2^{\f{n-r}2-1}
\\&=\f n{2i}\cdot
\f{(2x-1+i)^{n-1}-(2x-1-i)^{n-1}}{2^{n-1}}.\endalign$$ This yields
the result. \pro{Theorem 2.4} Let $m\in\Bbb N$ and
$n\in\{2,4,6,\ldots\}.$ Then
$$\aligned&\sum\Sb k=0\\4\mid k-2\endSb^n\b nk(-1)^{\f{k-2}4}2^{n-\f
k2}m^kB_{n-k}\\&=\f n2\Big\{m^{n-1}\Big((-1)^{[\f{n-2}4]}2^{\f
n2}-(1+(-1)^m)(-1)^{\f
n2}\Big)\\&\q+2i\sum_{r=1}^{[\f{m-1}2]}\Big((2r-m-mi)^{n-1}
-(2r-m+mi)^{n-1}\Big)\Big\}.\endaligned$$
\endpro
Proof. By Raabe's theorem, $\sum_{r=0}^{m-1}B_{n-k}\ls
rm=m^{k-(n-1)}B_{n-k}.$ Thus, using Theorem 2.3 we see that
$$\aligned&\sum\Sb k=0\\4\mid k-2\endSb^n\b nk(-1)^{\f{k-2}4}2^{n-\f
k2}m^kB_{n-k}\\&=\sum\Sb k=0\\4\mid k-2\endSb^n\b
nk(-1)^{\f{k-2}4}2^{n-\f k2}m^{n-1}\sum_{r=0}^{m-1}B_{n-k}\Ls
rm\\&=\f n2
i\sum_{r=0}^{m-1}\Big\{(2r-m-mi)^{n-1}-(2r-m+mi)^{n-1}\Big\}\\&=\f
n2i\Big\{(-m-mi)^{n-1}-(-m+mi)^{n-1}\\&\ \ +
\sum_{r=1}^{m-1}\big((2r-m-mi)^{n-1}-(2r-m+mi)^{n-1}\big)\Big\}\\&=\f
n2i\Big\{(-m-mi)^{n-1}-(-m+mi)^{n-1}-\big(1+(-1)^m\big)(mi)^{n-1}\\&\
\ +
2\sum_{r=1}^{[\f{m-1}2]}\big((2r-m-mi)^{n-1}-(2r-m+mi)^{n-1}\big)\Big\}\endaligned$$
It is evident that
$$\aligned &i\big((-1-i)^{n-1}-(-1+i)^{n-1}\big)\\&=i\big((-1-i)(2i)^{\f{n-2}2}
-(-1+i)(-2i)^{\f{n-2}2}\big)\\&=\cases
i\cdot(-2i)(2i)^{\f{n-2}2}=(-1)^{\f{n-2}4}2^{\f n2}&\t{if}\ 4\mid n-2,\\
i\cdot(-2)(2i)^{\f{n-2}2}=(-1)^{[\f{n-2}4]}2^{\f n2}&\t{if}\ 4\mid
n\endcases\endaligned$$ and
$$\align&(2(m-r)-m-mi)^{n-1}-(2(m-r)-m+mi)^{n-1}
\\&=(2r-m-mi)^{n-1}-(2r-m+mi)^{n-1}.\endalign$$
 Now combining all the above we
deduce the result.
\par Putting $m=1,2$ in Theorem 2.4 we deduce the following
identities.\pro{Corollary 2.4}(Ramanujan [R]) For $n=2,4,6,\ldots$
we have $$\sum\Sb k=0\\4\mid k-2\endSb^n\b
nk(-1)^{\f{k-2}4}2^{\f{n-k}2}B_{n-k}=(-1)^{[\f{n-2}4]}\f
n2.$$\endpro \pro{Corollary 2.5} For $n=2,4,6,\ldots$ we have
$$\sum\Sb k=0\\4\mid k-2\endSb^n\b
nk(-1)^{\f{k-2}4}2^{\f k2}B_{n-k}=\f
n2\big((-1)^{[\f{n-2}4]}2^{\f{n-2}2}+(-1)^{\f{n-2}2}\big).$$\endpro
\pro{Corollary 2.6} For $n=2,4,6,\ldots$ we have
$$\sum\Sb k=0\\4\mid k-2\endSb^n\b
nk(-1)^{\f{k-2}4}2^{\f{3k}2}B_{n-k}=\f
n2\big((-1)^{[\f{n-2}4]}2^{\f{3n}2-2}-(-1)^{\f
n2}2^{n-1}+4U_{n-1}(2,5)\big).$$ \endpro
Proof. Putting $m=4$ in
Theorem 2.4 we see that
$$\align&\sum\Sb k=0\\4\mid
k-2\endSb^n\b nk(-1)^{\f{k-2}4}2^{n+\f{3k}2}B_{n-k} \\&=\f
n2\Big\{4^{n-1}\big((-1)^{[\f{n-2}4]}2^{\f n2}-2(-1)^{\f
n2}\big)+2i\big((-2-4i)^{n-1}-(-2+4i)^{n-1}\big)\Big\}.\endalign$$
As
$$\aligned
&2i\big((-2-4i)^{n-1}-(-2+4i)^{n-1}\big)\\&=-2^ni\big((1+2i)^{n-1}-(1-2i)^{n-1}\big)\\&
=4\cdot 2^n\cdot\f 1{\sqrt{2^2-4\times
5}}\Big\{\Ls{2+\sqrt{2^2-4\times 5}}2^{n-1}-\Ls{2-\sqrt{2^2-4\times
5}}2^{n-1}\Big\}\\&=2^{n+2}U_{n-1}(2,5),\endaligned$$ the result
follows from the above.
\par Similarly, putting $m=3$ in Theorem 2.4 we deduce the
following result. \pro{Corollary 2.7} For $n=2,4,6,\ldots$ we have
$$\sum\Sb k=0\\4\mid k-2\endSb^n\b
nk(-1)^{\f{k-2}4}2^{n-\f k2}\cdot 3^kB_{n-k}=\f
n2\big((-1)^{[\f{n-2}4]}2^{\f n2}\cdot
3^{n-1}+12U_{n-1}(2,10)\big).$$ \endpro
 \pro{Theorem 2.5} For
$n\in\Bbb N$ we have
$$\align&\sum\Sb k=1\\4\mid k\endSb^n\b
nkB_{n-k}(x)((-1)^{\f k4}2^{\f k2-1}-1)\\&=\f
n8\Big\{2(x-1)^{n-1}-2x^{n-1}+(x+i)^{n-1}
+(x-i)^{n-1}\\&\qq-(x-1+i)^{n-1}-(x-1-i)^{n-1}\Big\}.\endalign$$
\endpro
Proof. From Lemma 2.1 we have
$$\align &\sum_{k=0}^{n-1}\b nkB_k(x)((1\pm i)^{n-k}-(\pm i)^{n-k}-(\pm i)^{n-k}
+(\pm i-1)^{n-k})
\\&=n(x\pm i)^{n-1}-n(x\pm i-1)^{n-1}.\endalign$$
Thus,
$$\align &2\sum\Sb k=0\\2\mid n-k\endSb^{n-1}\b
nkB_k(x)((i+1)^{n-k}-2i^{n-k}+(i-1)^{n-k})
\\&=\sum_{k=0}^{n-1}\b nkB_k(x)(1+(-1)^{n-k})((i+1)^{n-k}-2i^{n-k}+(i-1)^{n-k})
\\&=\sum_{k=0}^{n-1}\b nkB_k(x)((i+1)^{n-k}-2i^{n-k}
+(i-1)^{n-k}
\\&\qq+(-i+1)^{n-k}-2(-i)^{n-k}+(-i-1)^{n-k})
\\&=n(x+i)^{n-1}-n(x+i-1)^{n-1}+n(x-i)^{n-1}-n(x-i-1)^{n-1}.\endalign$$
As
$$\aligned (i+1)^{2m}-2i^{2m}+(i-1)^{2m}&=(2i)^m-2(-1)^m+(-2i)^m
\\&=\cases 2((-1)^{\f m2}2^m-1)&\t{if $2\mid m$,}
\\2&\t{if $2\nmid m$},\endcases\endaligned\tag 2.3$$
from the above we deduce
$$\align &4\sum\Sb k=0\\4\mid n-k\endSb^{n-1}\b
nkB_k(x)\big((-1)^{\f{n-k}4}2^{\f{n-k}2}-1\big)+4\sum\Sb k=0\\4\mid
n-k-2\endSb^{n-1}\b nkB_k(x)
\\&=n\big\{(x+i)^{n-1}+(x-i)^{n-1}-(x-1+i)^{n-1}-(x-1-i)^{n-1}\big\}.
\endalign$$
Using Lemma 2.1 we see that
$$\align &4\sum\Sb k=0\\4\mid
n-k\endSb^{n-1}\b nkB_k(x)+4\sum\Sb k=0\\4\mid n-k-2\endSb^{n-1}\b
nkB_k(x)\\&=4\sum\Sb k=0\\2\mid n-k\endSb^{n-1}\b nkB_k(x)
=2\sum_{k=0}^{n-1}\b nkB_k(x)(1^{n-k}-0^{n-k}-0^{n-k}+(-1)^{n-k})
\\&=2nx^{n-1}-2n(x-1)^{n-1}.\endalign$$
Therefore
$$\align &4\sum\Sb k=0\\4\mid n-k\endSb^{n-1}\b
nkB_k(x)\big((-1)^{\f{n-k}4}2^{\f{n-k}2}-2\big)+2nx^{n-1}-2n(x-1)^{n-1}
\\&=n\big\{(x+i)^{n-1}
+(x-i)^{n-1}-(x-1+i)^{n-1}-(x-1-i)^{n-1}\big\}.\endalign$$ To see
the result, we note that
$$\sum\Sb k=1\\4\mid k\endSb^n\b
nkB_{n-k}(x)\big((-1)^{\f k4}2^{\f k2}-2\big)=\sum\Sb k=0\\4\mid
n-k\endSb^{n-1}\b nkB_k(x)\big((-1)^{\f{n-k}4}2^{\f{n-k}2}-2\big).$$

 \pro{Corollary 2.8} For
$n=2,4,6,\ldots$ we have
$$\sum\Sb k=1\\4\mid k\endSb^n\b
nk((-4)^{\f k4}-2)B_{n-k}=\f n2\big((-1)^{[\f n4]}2^{\f
n2-1}-1\big).$$
\endpro
Proof. Note that $2\nmid n-1$ and
$$(-1+i)^{n-1}+(-1-i)^{n-1}=(-1+i)(-2i)^{\f{n-2}2}+(-1-i)(2i)^{\f{n-2}2}=-(-1)^{[\f
n4]}2^{\f n2}.$$ Now putting $x=0$ in Theorem 2.5 we deduce the
result.

 \subheading{3. Recurrence
formulas with gaps for Euler polynomials}
 \pro{Lemma 3.1} For any nonnegative integer $n$ we have
$$\sum_{k=0}^n\b
nkE_k(x)\big(z^{n-k}+(1+z)^{n-k}\big)=2(x+z)^n.$$
\endpro
Proof. It is well known that
$$E_n(x+y)=\sum_{k=0}^n\b nk E_n(x)y^{n-k}\qtq{and}
E_n(x)+E_n(x+1)=2x^n.$$ Thus,
$$\sum_{k=0}^n\b
nkE_k(x)\big(z^{n-k}+(1+z)^{n-k}\big)=E_n(x+z)+E_n(x+z+1)=2(x+z)^n.$$
This is the result. \pro{Theorem 3.1} For $n\in\Bbb N$ we have
$$\sum\Sb k=0\\4\mid k\endSb^n\b nk(-1)^{\f k4}2^{-\f
k2}E_{n-k}(x)=\f 12\Big\{\Big(x-\f 12+\f 12i\Big)^n+\Big(x-\f 12-\f
12i\Big)^n\Big\}.$$
\endpro Proof. Putting $z=\f 12(\pm i-1)$ in
Lemma 3.1 we obtain $$\sum_{k=0}^n\b nkE_k(x)\ls 12^{n-k}\big((\pm
i-1)^{n-k}+(1\pm i)^{n-k}\big)=2\Big(x+\f 12(\pm i-1)\Big)^n.$$
 Thus $$\aligned &\sum_{k=0}^n\b nkE_k(x)\ls
12^{n-k}\\&\qq\times\Big\{(i-1)^{n-k}+(1+i)^{n-k}+(-i-1)^{n-k}+(1-i)^{n-k}\Big\}
\\&=2\Big(x-\f 12+\f 12i\Big)^n+2\Big(x-\f 12-\f 12i\Big)^n.\endaligned$$ Clearly
$$\aligned
&(i-1)^{n-k}+(1+i)^{n-k}+(-i-1)^{n-k}+(1-i)^{n-k}\\&=\big(1+(-1)^{n-k}\big)\big((i-1
)^{n-k}+(1+i)^{n-k}\big)\\&=\big(1+(-1)^{n-k}\big)\big((-2i
)^{\f{n-k}2}+(2i)^{\f{n-k}2}\big)\\&=\cases
2^{\f{n-k}2+2}(-1)^{\f{n-k}4}&\t{if} \ 4\mid n-k,\\0&\t{if} \ 4\mid
n-k-2.\endcases\endaligned$$ So $$\aligned&\sum\Sb k=0\\4\mid
k\endSb^n\b nkE_{n-k}(x)\ls 12^k(-1)^{\f k4}2^{\f k2+2}\\&=\sum\Sb
k=0\\4\mid n-k\endSb^n\b nkE_k(x)\ls
12^{n-k}(-1)^{\f{n-k}4}2^{\f{n-k}2+2}\\&=2\Big\{\Big(x-\f 12+\f
12i\Big)^n+\Big(x-\f 12-\f 12i\Big)^n\Big\}.\endaligned$$ This
yields the result. \pro{Corollary 3.1 ([L, (15)])} For
$n=2,4,6,\ldots$ we have
$$\sum\Sb k=0\\4\mid k\endSb^n\b nk(-4)^{\f k4}E_{n-k}=(-1)^{\f
n2}.$$
\endpro
Proof. As $E_n=2^nE_n(\f 12)$, taking $x=\f 12$ in Theorem 3.1 we
deduce the result.
 \pro{Corollary 3.2} For $n=2,4,6,\ldots$ we have
$$\sum\Sb k=0\\4\mid k\endSb^{n-1}\b nk(-1)^{\f
k4}2^{\f {n-k}2}(2^{n-k}-1)B_{n-k}=(-1)^{[\f n4]}\f n2.$$
\endpro
Proof. Taking $x=0$ and replacing $n$ with $n-1$ in Theorem 3.1 and
then applying (1.4) we deduce the result.\par\q
\par We remark that Corollary 3.2 is equivalent to [S, Corollary
4.1].
 \pro{Theorem 3.2} For any nonnegative integer $n$
we have
$$\align &\sum\Sb k=0\\4\mid k\endSb^n\b nk
\big((-1)^{\f k4}2^{\f k2-1}+1\big)E_{n-k}(x)
\\&=\f
14\big(2x^n+2(x-1)^n+(x+i)^n+(x-i)^n+(x-1+i)^n+(x-1-i)^n\big).\endalign$$
\endpro

Proof. From Lemma 3.1 we have
$$\sum_{k=0}^n\b nkE_k(x)((1\pm i)^{n-k}+2(\pm
i)^{n-k}+(-1\pm i)^{n-k})=2(x\pm i)^n+2(x-1\pm i)^n.$$ Thus,
$$\align &2\sum\Sb k=0\\2\mid n-k\endSb^n\b
nkE_k(x)((i+1)^{n-k}+2i^{n-k}+(i-1)^{n-k})
\\&=\sum_{k=0}^n\b nkE_k(x)(1+(-1)^{n-k})((i+1)^{n-k}+2i^{n-k}+(i-1)^{n-k})
\\&=\sum_{k=0}^n\b nkE_k(x)((i+1)^{n-k}+2i^{n-k}
+(i-1)^{n-k}
\\&\qq+(-i+1)^{n-k}+2(-i)^{n-k}+(-i-1)^{n-k})
\\&=2(x+i)^n+2(x+i-1)^n+2(x-i)^n+2(x-i-1)^n.\endalign$$ As
$$\aligned (i+1)^{2m}+2i^{2m}+(i-1)^{2m}&=(2i)^m+2(-1)^m+(-2i)^m
\\&=\cases 2((-1)^{\f m2}2^m+1)&\t{if $2\mid m$,}
\\-2&\t{if $2\mid m-1$},\endcases\endaligned$$
from the above we deduce
$$\align &4\sum\Sb k=0\\4\mid n-k\endSb^n\b
nkE_k(x)\big((-1)^{\f{n-k}4}2^{\f{n-k}2}+1\big)-4\sum\Sb k=0\\4\mid
n-k-2\endSb^n\b nkE_k(x)
\\&=2\big\{(x+i)^n+(x-1+i)^n+(x-i)^n+(x-1-i)^n\big\}.
\endalign$$
Using Lemma 3.1 we see that
$$\align &4\sum\Sb k=0\\4\mid
n-k\endSb^n\b nkE_k(x)+4\sum\Sb k=0\\4\mid n-k-2\endSb^n\b
nkE_k(x)\\&=4\sum\Sb k=0\\2\mid n-k\endSb^n\b nkE_k(x)
=2\sum_{k=0}^n\b
nkE_k(x)(1^{n-k}+0^{n-k}+0^{n-k}+(-1)^{n-k})
\\&=2\big(E_n(x+1)+E_n(x)+E_n(x)+E_n(x-1)\big)
\\&=4x^n+4(x-1)^n.\endalign$$
Therefore
$$\align &4\sum\Sb k=0\\4\mid n-k\endSb^n\b
nkE_k(x)\big((-1)^{\f{n-k}4}2^{\f{n-k}2}+2\big)
\\&=2\big((x+i)^n+(x+i-1)^n+(x-i)^n+(x-i-1)^n+2x^n+2(x-1)^n
\big).\endalign$$ Replacing $k$ with $n-k$ in the above formula we
derive the result.

 \pro{Theorem 3.3} For
any nonnegative integer $n$ we have
$$\align&4E_n(x)+3\sum_{k=1}^{[n/6]} \b
n{6k}E_{n-6k}(x)\\&=x^n+(x-1)^n+(-1)^nV_n(1-2x,x^2-x+1).\endalign$$
\endpro Proof.  Taking
$z=\o,\o^2$ in Lemma 3.1 we see that
$$\sum_{k=0}^n\b nkE_k(x)\big(\o^{n-k}+(1+\o)^{n-k}\big)
=2(x+\o)^n$$ and
$$\sum_{k=0}^n\b nkE_k(x)\big(\o^{2(n-k)}+(1+\o^2)^{n-k}\big)
=2(x+\o^2)^n.$$ Hence

$$\align
&\sum_{k=0}^n\b
nkE_k(x)\big(\omega^{n-k}+{\omega}^{2(n-k)}+(-\omega^2)^{n-k}
+(-\omega)^{n-k}\big)
\\&=2(x+w)^n+2(x+\o^2)^n.\endalign$$ Since
$$\aligned&\omega^{n-k}+{\omega}^{2(n-k)}+(-\omega^2)^{n-k}
+(-\omega)^{n-k}\\&=(1+(-1)^{n-k})(\omega^{n-k}+\omega^{2(n-k)})=\cases
0&\t{if $2\nmid n-k$,}
\\ 4&\t{if $6\mid n-k$,}
\\2(\omega+\omega^2)=-2&\t{if $ 6\mid n-k\pm 2$,}\endcases\endaligned$$
we see that $$\align &6\sum\Sb k=0\\6\mid n-k\endSb^n\b nkE_k(x)
-2\sum\Sb k=0\\2\mid n-k\endSb^n\b nkE_k(x)
\\&=\sum_{k=0}^n\b nkE_k(x)(1+(-1)^{n-k})(\omega^{n-k}+\omega^{2(n-k)})
\\&=2(x+\omega)^n+2(x+\omega^2)^n.\endalign$$
That is,
$$6\sum\Sb k=0\\6\mid k\endSb^n\b nkE_{n-k}(x)
-2\sum\Sb k=0\\2\mid k\endSb^n\b nkE_{n-k}(x)
=2(x+\omega)^n+2(x+\omega^2)^n.$$ By Lemma 3.1 we have
$$\align 2E_n(x)+2\sum\Sb k=0\\2\mid n-k\endSb^n\b nkE_k(x)
&=\sum_{k=0}^n\b
nkE_k(x)(0^{n-k}+1^{n-k}+(-1)^{n-k}+0^{n-k})\\&=2x^n+2(x-1)^n.\endalign$$Thus,
$$\sum\Sb k=0\\2\mid k\endSb^n\b nkE_{n-k}(x)
=\sum\Sb k=0\\2\mid n-k\endSb^n\b nkE_k(x)=x^n+(x-1)^n-E_n(x).$$
Hence,
$$\align &6\sum_{k=0}^{[n/6]}\b n{6k}E_{n-6k}(x)
\\&=2\sum_{k=0}^{[n/2]}\b n{2k}E_{n-2k}(x)+2(x+\omega)^n+2(x+\omega^2)^n
\\&=2x^n+2(x-1)^n-2E_n(x)+2(x+\omega)^n+2(x+\omega^2)^n.
\endalign$$
Therefore,
$$4E_n(x)=x^n+(x-1)^n+(x+\omega)^n+(x+\omega^2)^n-3\sum_{k=1}^{[n/6]}
\b n{6k}E_{n-6k}(x).$$ To see the result, we note that
$$\align (-1)^n\big((x+\omega)^n+(x+\omega^2)^n
\big)&=\Ls{1-2x-\sqrt{-3}}2^n+\Ls{1-2x+\sqrt{-3}}2^n
\\&=\Ls{1-2x+\sqrt{(1-2x)^2-4(x^2-x+1)}}2^n
\\&\q+\Ls{1-2x-\sqrt{(1-2x)^2-4(x^2-x+1)}}2^n
\\&=V_n(1-2x,x^2-x+1).
\endalign$$
\par \q\par Note that $E_n=2^nE_n(\f 12)$. Putting $x=\f 12$ in Theorem 3.3
we deduce (1.6).

\enddocument

\Refs\widestnumber\key{BEW}
 \ref\key C\by M. Chellali
\paper Acc$\acute e$l$\acute e$ration de calcul de nombres de
Bernoulli\jour J. Number Theory\vol 28\yr 1988\pages 347-362\endref

 \ref\key L\by D.H. Lehmer\paper Lacunary recurrence formulas for
the numbers of Bernoulli and Euler\jour Ann. Math. \vol 36\yr
1935\pages 637-649\endref

\ref \key MOS\by W. Magnus, F. Oberhettinger and R.P. Soni\book
Formulas and Theorems for the Special Functions of Mathematical
Physics $(3rd\ Edition)$\publ Springer-Verlag\publaddr New York\yr
1966\pages 25-32\endref

 \ref\key S \by Z.H. Sun\paper On the properties of
Newton-Euler pairs\jour J. Number Theory\vol 114\yr 2005\pages
88-123\endref

\ref\key R\by S. Ramanujan\paper Some properties of Bernoulli's
numbers \jour J. Indian Math. Soc.\vol 3\yr 1911\pages
219-234\endref \ref\key W\by
 H. C. Williams\book \'Edouard Lucas and Primality
 Testing, Canadian Mathematical Society Series of Monographs and
 Advanced Texts (Vol.22) \publ Wiley, New York \yr 1998\pages
  74-92\endref


\endRefs
\bye